\documentclass[11pt,reqno]{amsart} 
\usepackage{amsmath, amsthm, amsfonts, amssymb} 
\usepackage[colorlinks]{hyperref} 
\usepackage[T1]{fontenc} 
\usepackage[mathscr]{eucal} 
\usepackage[polish,english]{babel}
\usepackage{graphicx}
\usepackage{float}
\newtheorem{theorem}{Theorem}[section]

\newtheorem{definition}[theorem]{Definition}

\begin{document} \title[Exploring the EHM Index of Graphs: Application and Prediction] {Exploring the Edge Hyper-Zagreb Index of Graphs: Applications and Prediction of Thermodynamic Properties for Organic Linear Acenes Molecules} 
\date{} \author [Z. Aliannejadi and S. Shafiee Alamoti]{Z. Aliannejadi and S. Shafiee Alamoti} 
\begin{abstract}
The paper discusses the edge hyper-Zagreb index of a graph, which is calculated by replacing vertex degrees with edge degrees. The degree of an edge is determined by adding up the degrees of the end vertices of the edge and subtracting $2$. We examine the edge hyper-Zagreb index of the Cartesian product and join of graphs, and also calculate it for organic linear Acenes molecules with the formula $(C_{4n+2}H_{2n+4})$. We establish a correlation between topological indices based on the number of rings and predict thermodynamic properties of Acenes family, such as electron affinity, bond, heat of formation and gap energy, using the Topological Indices Method $(TIM)$. 

\end{abstract} 
\maketitle \textbf{Mathematics Subject Classification (2010): 05C35; 05C07.}
\section{Introduction} 
A topological index, alternatively referred to as a connectivity index, is a form of molecular descriptor employed in chemical graph theory, molecular topology, and mathematical chemistry. It is computed from the molecular graph of a chemical compound and serves as a numerical parameter that describes the graph's topology [17]. Typically, topological indices remain invariant across graphs and are integral in the construction of quantitative structure-activity relationships $(QSARs)$ that establish connections between the biological activity or various properties of molecules and their chemical structure [16].\\
To compute topological indices, hydrogen-depleted molecular graphs are employed, where atoms are depicted as vertices and bonds as edges. The interconnections among the atoms are delineated using diverse types of topological matrices, such as distance or adjacency matrices, which are manipulated mathematically to yield a single numerical value identified as a graph invariant, graph-theoretical index, or topological index [13,14]. These indices, which are $2D$ descriptors, can be easily computed from the molecular graphs and are independent of the graph's depiction or labeling. Furthermore, the structural energy minimization of the chemical compound is unnecessary. Topological indices find applications in chemical graph theory, molecular topology, and mathematical chemistry for delineating the molecular graph's topology.\\ 
Basic topological indices do not consider double bonds or atom characteristics, and they disregard hydrogen atoms [22]. They apply exclusively to connected undirected molecular graphs. Advanced topological indices consider the hybridization status of individual atoms within the molecule [15]. The Hosoya index is the initial and most prevalent topological index, while others encompass the Wiener index, Randi{\'c}'s molecular connectivity index, Balaban's J index [7], and the $TAU$ descriptors [25]. These indices may yield identical values for distinct molecular graphs, hence utilizing multiple indices can enhance the discriminatory capacity.\\
The enhanced topochemical atom $(ETA)$ [23] indices were created through the enhancement of $TAU$ descriptors. A topological index may yield the same result for a subset of different molecular graphs, meaning it cannot distinguish between graphs within this subset. The ability to differentiate between graphs is a crucial characteristic of a topological index. To improve this capability, multiple topological indices can be merged to form a comprehensive superindex [9]. $QSARs$ are models used for prediction that are developed by using statistical methods to correlate the biological activity (both positive therapeutic effects and negative side effects) of substances (such as drugs, toxicants, and environmental pollutants) with descriptors that represent molecular structure and/or properties. $QSARs$ find applications in various fields including risk assessment, toxicity prediction, regulatory decision-making [27], as well as in drug discovery and lead optimization [21].\\
Nanotechnology encompasses the manipulation of matter at the atomic, molecular, and supramolecular scales [28]. Initially, the goal of nanotechnology was to manipulate atoms and molecules precisely in order to create larger products, a concept known as molecular nanotechnology [1,2,3,4,18]. Through nanotechnology, the production, manipulation, and utilization of nanomaterials are made possible [17]. Recently, there has been growing interest in organic molecules like Acenes, which exhibit nanostructure dimensions [3,7,14,25,26]. Investigating the physical, chemical, electrical, and thermodynamic properties of heavier Acenes is a challenging and time-intensive endeavor. Establishing a reliable model for predicting these properties is essential. One effective approach to accomplishing this is by generating a topological index for the chemical molecules [24].\\
In this paper, we focus on a connected, finite and undirected graphs. The degree of a vertex $a$ in a graph $\Gamma$, denoted by $d_\Gamma (a)$, is the number of edges that connect it.\\ 
Gutman and Trinajstic [19] introduced the first and second kind of Zagreb indices in $1972$, which were used to study the structure dependency of the total $varphi$-electron energy. The first and second Zagreb indices of a connected graph $\Gamma$ for $\alpha=ab\in E(\Gamma)$ denoted by $EM_1(\Gamma)$ and $EM_2(\Gamma)$, are defined as
\begin{center}
$M_1(\Gamma)=\sum_{a\in V(\Gamma)}{(d_\Gamma (a))}^2=\sum_{ab \in E(\Gamma)}(d_\Gamma (a)+d_\Gamma (b))$
\end{center}
and
\begin{center}
$M_2(\Gamma)=\sum _{ab \in E(\Gamma)}(d_\Gamma (a).d_\Gamma (b))$,
\end{center}
respectively. Ashrafi et al. [6] introduced the first and second Zagreb coindices, 
\begin{center}
$\overline{M}_1(\Gamma)=\sum_{ab \notin E(\Gamma)}(d_\Gamma (a)+d_\Gamma (b))$ 
\end{center}
and
\begin{center}
$\overline{M}_2(\Gamma)=\sum _{ab \notin E(\Gamma)}(d_\Gamma (a).d_\Gamma (b))$,
\end{center}
respectively. These indices were found to be related to the total $\varphi$-electron energy [19] and were proposed as measures of branching of the carbon-atom skeleton [20]. The Zagreb index, or Zagreb group index, was later named in a review article [8]. Milicevic et al. [5] reformulated the Zagreb indices in terms of edge-degrees instead of vertex-degrees, as below
\begin{center}
$EM_1(\Gamma)=\sum_{\alpha\in E(\Gamma)}{(d_\Gamma (\alpha))}^2$ $\hspace*{1cm}$ $EM_2(\Gamma)=\sum_{\alpha\thicksim \beta}(d_\Gamma (\alpha)).(d_\Gamma (\beta))$,
\end{center}
where $d_\Gamma (\alpha)$ denotes the degree of the edge $\alpha$ in $\Gamma$, which is defined by 
$d_\Gamma (\alpha)=d_\Gamma (a)+d_\Gamma (b)-2$ with $\alpha=ab$ and $\alpha\thicksim \beta$ means that the edges $\alpha$ and $\beta$ are adjacent, i.e., they share a common end-vertex in $\Gamma$.\\
Shirdel et al. [11] proposed a new topological index called the hyper Zagreb index, which is defined for a connected $\Gamma$ as
\begin{center}
$HM(\Gamma)=\sum_{ab \in E(\Gamma)}(d_\Gamma (a)+d_\Gamma (b))^2 $.\\
\end{center}
Now, We would like to define the edge hyper-Zagreb index as:
\begin{definition}
A edge hyper-Zagreb index of graph $\Gamma$ for $\alpha=ab\in E(\Gamma)$ is introduced as
\begin{center}
$EHM(\Gamma)=\sum_{\alpha\thicksim \beta\in E(\Gamma)}{(d_\Gamma (\alpha)+d_\Gamma (\alpha))}^2$.
\end{center}
\end{definition}
\section{Cartesian Product and Join of Graphs}
In the following section, we explore the behavior of the edge hyper-Zagreb index under the Cartesian product and join of graphs. We consider only connected, finite and simple graphs, denoted by $\Gamma$ and $\Omega$ with vertex sets $V(\Gamma)$ and $V(\Omega)$ and edge sets $E(\Gamma)$ and $E(\Omega)$, respectively. In addition, we symbolize $|V(\Gamma)|=n_\Gamma$, $|E(\Gamma)|=m_\Gamma$ and $|V(\Omega)|=n_\Omega$, $|E(\Omega)|=m_\Omega$.\\\\
The join of $\Gamma +\Omega$ of graphs $\Gamma$ and $\Omega$ is a graph with vertex set $V(\Gamma) \bigcup V(\Omega)$ and edge set $E(\Gamma) \bigcup E(\Omega) \bigcup \{ab|a\in V(\Gamma)$ and $b\in V(\Omega)\}$. The degree of a vertex $a$ of $\Gamma +\Omega$ is given by\\
\begin{equation} d_{(\Gamma +\Omega)}(a)=
\left\{
\begin{array}{cc}
d_{\Gamma}(a)+n_\Omega, a\in V(\Gamma) \notag\\
d_{\Omega}(a)+n_\Gamma, a\in V(\Omega).\\
\end{array}
\right.
\end{equation}
The Cartesian product $\Gamma \times \Omega$ of graphs $\Gamma$ and $\Omega$ has the vertex set $V(\Gamma \times \Omega)=V(\Gamma) \times V(\Omega)$ and $(a,b)(\acute{a},\acute{b})$ is an edge of $\Gamma \times \Omega$ if $a=\acute{a}$ and $b\acute{b}\in E(\Omega)$, or $a\acute{a}\in E(\Gamma)$ and $b=\acute{b}$.\\
\begin{theorem}
The edge hyper-Zagreb index of $\Gamma +\Omega$ is equal to:\\

$EHM(\Gamma +\Omega)= 2EM_1(\Gamma)+16n_\Omega(M_1(\Gamma)-2m_\Gamma)+16n_\Omega^2m_\Gamma$\\
$\hspace*{2.75cm}$ $+2EM_2(\Gamma) + 2EM_1(\Omega)+16n_\Gamma (M_1(\Omega)-2m_\Omega)$\\
$\hspace*{2.75cm}$ $+16n_\Gamma^2m_\Omega+2EM_2(\Omega) +
2n_\Omega M_1(\Gamma) + 2n_\Gamma M_1(\Omega)$\\
$\hspace*{2.75cm}$ $+ 32m_\Gamma m_\Omega + 8(2n_\Gamma + 2n_\Omega - 4)(n_\Gamma m_\Omega + n_\Omega m_\Gamma)$\\ 
$\hspace*{2.75cm}$ $+ 4n_\Gamma n_\Omega (2n_\Gamma + 2n_\Omega - 4)^2 + 8{m_\Gamma}^2 + 8{m_\Omega}^2$.\\

\end{theorem}

\begin{proof}
Let $\alpha=ab$ and $\beta=bc$ be two adjacent edges in $E(\Gamma +\Omega)$. Then we have\\

$EHM(\Gamma +\Omega)=\sum_{\alpha\thicksim\beta \in E(\Gamma +\Omega)}(d_{\Gamma +\Omega}(\alpha) + d_{\Gamma +\Omega}(\beta))^2$\\

$\hspace*{2.5cm}$=$\sum_{\alpha\thicksim\beta \in E(\Gamma)}(d_{\Gamma +\Omega}(\alpha) + d_{\Gamma +\Omega}(\beta))^2$\\

$\hspace*{2.4cm}$ $+\sum_{\alpha\thicksim\beta \in E(\Omega)}(d_{\Gamma +\Omega}(\alpha) + d_{\Gamma +\Omega}(\beta))^2$\\

$\hspace*{2.4cm}$ $+\sum_{\alpha\thicksim\beta \in \{\alpha\thicksim\beta|\alpha\in E(\Gamma), \beta\in E(\Omega)\}}(d_{\Gamma +\Omega}(\alpha) + d_{\Gamma +\Omega}(\beta))^2.$\\

$\hspace*{2.4cm}$ $=A1+A_2+A_3$.

We consider\\

$\hspace*{1.8cm}$ $A_1=\sum_{\alpha\thicksim\beta \in E(\Gamma)}(d_{\Gamma +\Omega}(\alpha) + d_{\Gamma +\Omega}(\beta))^2$\\

$\hspace*{2.4cm}$=$\sum_{\alpha\thicksim\beta \in E(\Gamma)}(d_{\Gamma +\Omega}(a) +d_{\Gamma +\Omega}(b) -2 + d_{\Gamma +\Omega}(b) + d_{\Gamma +\Omega}(c) - 2)^2$\\

$\hspace*{2.28cm}$=$\sum_{\alpha\thicksim\beta \in E(\Gamma)}(d_{\Gamma}(a) +d_{\Gamma}(b) +2n_\Omega - 2 + d_{\Gamma}(b) + d_{\Gamma}(c) + 2n_\Omega - 2)^2$\\

$\hspace*{2.4cm}$=$\sum_{\alpha\thicksim\beta \in E(\Gamma)}(d_{\Gamma}(a) +d_{\Gamma}(b) - 2 + d_{\Gamma}(b) + d_{\Gamma}(c) - 2 + 4n_\Omega)^2$\\

$\hspace*{2.4cm}$=$\sum_{\alpha\thicksim\beta \in E(\Gamma)}((d_{\Gamma}(a) +d_{\Gamma}(b) -2)^2
+ (d_{\Gamma}(b) + d_{\Gamma}(c) - 2)^2  $\\

$\hspace*{2.3cm}$ $+ 16n_\Omega^2 + 2(d_{\Gamma}(a) +d_{\Gamma}(b) -2)(d_{\Gamma}(b) + d_{\Gamma}(c) - 2)$\\

$\hspace*{2.3cm}$ $+8n_\Omega (d_{\Gamma}(a) + d_{\Gamma}(b) - 2) + 8n_\Omega (d_{\Gamma}(b) + d_{\Gamma}(c) - 2))$\\

$\hspace*{2.3cm}$ $=\sum_{\alpha\thicksim\beta \in E(\Gamma)}(d_{\Gamma}(a) +d_{\Gamma}(b) -2)^2$\\

$\hspace*{2.3cm}$ $+ 8n_\Omega \sum_{\alpha\thicksim\beta \in E(\Gamma)} (d_{\Gamma}(a) +d_{\Gamma}(b) -2)$\\

$\hspace*{2.3cm}$ $+16\sum_{\alpha\thicksim\beta \in E(\Gamma)}n_\Omega^2 + \sum_{\alpha\thicksim\beta \in E(\Gamma)}(d_{\Gamma}(b) +d_{\Gamma}(c) -2)^2$\\

$\hspace*{2.3cm}$ $+8n_\Omega \sum_{\alpha\thicksim\beta \in E(\Gamma)}(d_{\Gamma}(b) +d_{\Gamma}(c) -2)^2$ \\

$\hspace*{2.3cm}$ $+2\sum_{\alpha\thicksim\beta \in E(\Gamma)}(d_{\Gamma}(a) +d_{\Gamma}(b) -2)(d_{\Gamma}(b) +d_{\Gamma}(c) -2)$\\

$\hspace*{2.3cm}$ $=EM_1(\Gamma)+8n_\Omega(M_1(\Gamma) -  2m_\Gamma) +16 n_\Omega^2 m_\Gamma$\\

$\hspace*{2.3cm}$ $+EM_1(\Gamma)+8n_\Omega(M_1(\Gamma) +2EM_2(\Gamma)$\\

$\hspace*{2.3cm}$ $=2EM_1(\Gamma)+16n_\Omega(M_1(\Gamma)-2m_\Gamma)+16n_\Omega^2m_\Gamma+2EM_2(\Gamma).$\\

Similarly,

$\hspace*{1.8cm}$ $A_2=\sum_{\alpha\thicksim\beta \in E(\Omega)}(d_{\Gamma +\Omega}(\alpha) + d_{\Gamma +\Omega}(\beta))^2$\\

$ \hspace*{2.1cm}$ $=\sum_{\alpha\thicksim\beta \in E(\Omega)}(d_{\Gamma +\Omega}(a) +d_{\Gamma +\Omega}(b) -2 + d_{\Gamma +\Omega}(b) + d_{\Gamma +\Omega}(c) - 2)^2$\\

$\hspace*{2cm}$ $=\sum_{\alpha\thicksim\beta \in E(\Omega)}(d_{\Omega}(a) +d_{\Omega}(b) +2n_\Gamma - 2 + d_{\Omega}(b) + d_{\Omega}(c) + 2n_\Gamma - 2)^2$\\

$\hspace*{2.1cm}$ $=\sum_{\alpha\thicksim\beta \in E(\Omega)}(d_{\Omega}(a) +d_{\Omega}(b) - 2 + d_{\Omega}(b) + d_{\Omega}(c) - 2 + 4n_\Gamma)^2$\\

$\hspace*{2.12cm}$ $=\sum_{\alpha\thicksim\beta \in E(\Omega)}((d_{\Omega}(a) +d_{\Omega}(b) -2)^2
+ (d_{\Omega}(b) + d_{\Omega}(c) - 2)^2 + 16n_\Gamma^2 $\\

$\hspace*{2.12cm}$ $+ 2(d_{\Omega}(a) +d_{\Omega}(b) -2)(d_{\Omega}(b) + d_{\Omega}(c) - 2)$\\

$\hspace*{2.12cm}$ $ +8n_\Gamma (d_{\Omega}(a) + d_{\Omega}(b) - 2) + 8n_\Gamma (d_{\Omega}(b) + d_{\Omega}(c) - 2))$\\

$\hspace*{2.12cm}$ $=\sum_{\alpha\thicksim\beta \in E(\Omega)}(d_{\Omega}(a) +d_{\Omega}(b) -2)^2$\\

$\hspace*{2.12cm}$ $+ 8n_\Gamma \sum_{\alpha\thicksim\beta \in E(\Omega)} (d_{\Omega}(a) +d_{\Omega}(b) -2)$\\

$\hspace*{2.1cm}$ $+16\sum_{\alpha\thicksim\beta \in E(\Omega)}n_\Gamma^2 + \sum_{\alpha\thicksim\beta \in E(\Omega)}(d_{\Omega}(b) +d_{\Omega}(c) -2)^2$\\

$\hspace*{2.1cm}$ $+8n_\Gamma \sum_{\alpha\thicksim\beta \in E(\Omega)}(d_{\Omega}(b) +d_{\Omega}(c) -2)^2$ \\

$\hspace*{2.1cm}$ $+2\sum_{\alpha\thicksim\beta \in E(\Omega)}(d_{\Omega}(a) +d_{\Omega}(b) -2)(d_{\Omega}(b) +d_{\Omega}(c) -2)$\\

$\hspace*{2.1cm}$ $=EM_1(\Omega)+8n_\Gamma (M_1(\Omega) -  2m_\Omega) +16 n_\Gamma^2 m_\Omega$\\

$\hspace*{2.1cm}$ $+EM_1(\Omega)+8n_\Gamma (M_1(\Omega) +2EM_2(\Omega)$\\

$\hspace*{2.1cm}$ $=2EM_1(\Omega)+16n_\Gamma (M_1(\Omega)-2m_\Omega)+16n_\Gamma^2m_\Omega+2EM_2(\Omega).$\\

Finally,

$\hspace*{1.5cm}$ $A_3=\sum_{\alpha\thicksim\beta \in \{\alpha\thicksim\beta|\alpha\in E(\Gamma), \beta\in E(\Omega)\}}(d_{\Gamma +\Omega}(\alpha) + d_{\Gamma +\Omega}(\beta))^2$

\begin{center}
$=\sum_{\alpha\thicksim\beta \in \{\alpha\thicksim\beta|\alpha\in E(\Gamma), \beta\in E(\Omega)\}}(d_{\Gamma +\Omega}(a) + d_{\Gamma +\Omega}(b) - 2 +d_{\Gamma +\Omega}(b) + d_{\Gamma +\Omega}(c) - 2  )^2$\\
\end{center}

\begin{center}
$=\sum_{a,b\in V(\Gamma) , b,c\in V(\Omega)}(d_{\Gamma}(a) + n_\Omega + d_{\Omega}(b) + n_\Gamma - 2 +d_{\Gamma}(b) + n_\Omega + d_{\Omega}(c) + n_\Gamma - 2  )^2$
\end{center}

\begin{center}
$=\sum_{a,b\in V(\Gamma) , b,c\in V(\Omega)}((d_{\Gamma}(a) + d_{\Omega}(b)) +(d_{\Gamma}(b) + d_{\Omega}(c)) + (2n_\Omega+ 2n_\Gamma - 4)^2$
\end{center}

\begin{center}
$=\sum_{a,b\in V(\Gamma) , b,c\in V(\Omega)}((d_{\Gamma}(a)^2 + d_{\Omega}(b)^2 +2 d_{\Gamma}(a)d_{\Omega}(b)) + (d_{\Gamma}(b)^2 + d_{\Omega}(c)^2+2 d_{\Gamma}(b) d_{\Omega}(b)) + (2n_\Omega + 2n_\Gamma - 4)^2$
\end{center}

\begin{center}
$+2 (d_{\Gamma}(a) + d_{\Omega}(b))(d_{\Gamma}(b) + d_{\Omega}(c)) $
\end{center}

\begin{center}
$+2 (d_{\Gamma}(a) + d_{\Omega}(b))(2n_\Omega + 2n_\Gamma- 4) $
\end{center}

\begin{center}
$+2 (d_{\Gamma}(b) + d_{\Omega}(a))(2n_\Omega + 2n_\Gamma- 4) $
\end{center}

$=n_\Omega M_1(\Gamma) + n_\Gamma M_1(\Omega) + 8m_\Gamma m_\Omega + n_\Omega M_1(\Gamma) + n_\Gamma M_1(\Omega)+ 8m_\Gamma m_\Omega $\\

 $+ 4n_\Gamma n_\Omega ( 2n_\Gamma + 2n_\Omega - 4)^2 + 4 (2n_\Gamma + 2n_\Omega - 4)( n_\Gamma m_\Omega + n_\Omega m_\Gamma) + 4 (2n_\Gamma$\\ 
 
 $+ 2n_\Omega - 4)( n_\Gamma m_\Omega + n_\Omega m_\Gamma) + 8m_\Gamma m_\Omega + 8m_\Gamma m_\Omega + 8{m_\Gamma} ^2 + 8{m_\Omega}^2$\\

$=2n_\Omega M_1(\Gamma) + 2n_\Gamma M_1(\Omega) + 32m_\Gamma m_\Omega + 8(2n_\Gamma + 2n_\Omega - 4)(n_\Gamma m_\Omega + n_\Omega m_\Gamma)$\\ 

$+ 4n_\Gamma n_\Omega (2n_\Gamma + 2n_\Omega - 4)^2 + 8{m_\Gamma}^2 + 8{m_\Omega}^2$.\\

Adding $A_1, A_2$ and $A_3$, we obtain the desired result.
\end{proof}

\begin{theorem}
The edge hyper-Zagreb index of $\Gamma \times \Omega$ is given by\\

$EHM(\Gamma \times \Omega)=2(n_\Gamma EM_1(\Omega)+8m_\Gamma M_1(\Omega)+4m_\Omega M_1(\Gamma)$\\
$\hspace*{2.8cm}$ $+ n_\Omega EM_1(\Gamma)+8m_\Omega M_1(\Gamma)+4m_\Gamma M_1(\Omega) - 32m_\Gamma m_\Omega)$.
\end{theorem}
\begin{proof}
Let $(a,b)(\acute{a},\acute{b})$ and $(\acute{a},\acute{b})(\acute{\acute{a}},\acute{\acute{b}})$ be two adjacent edges in $EHM(\Gamma \times \Omega)$. Then we have:\\

$EHM(\Gamma \times \Omega)=\sum_{(a,b)(\acute{a},\acute{b}) \thicksim (\acute{a},\acute{b})(\acute{\acute{a}},\acute{\acute{b}}) \in E(\Gamma \times \Omega)}(d_{\Gamma \times \Omega}(a,b)(\acute{a},\acute{b})+d_{\Gamma \times \Omega}(\acute{a},\acute{b})(\acute{\acute{a}},\acute{\acute{b}}))^2$ \\

$\hspace*{2.2cm}$ $=\sum_{(a,b)(\acute{a},\acute{b}) \thicksim (\acute{a},\acute{b})(\acute{\acute{a}},\acute{\acute{b}}) \in E(\Gamma \times \Omega)}(d_{\Gamma \times \Omega}(a,b)+d_{\Gamma \times \Omega}(\acute{a},\acute{b}) - 2 $\\

$\hspace*{2.2cm}$ $+ d_{\Gamma \times \Omega}(\acute{a},\acute{b}) + d_{\Gamma \times \Omega}(\acute{\acute{a}},\acute{\acute{b}}) - 2)^2$\\

$\hspace*{2.2cm}$ $=\sum_{(a,b)(a,\acute{b}) ,  b\acute{b} \in E(\Omega)}(d_{\Gamma \times \Omega}(a,b)+d_{\Gamma \times \Omega}(a,\acute{b}) - 2)^2 $\\

$\hspace*{2.2cm}$ $+\sum_{(a,b)(\acute{a},b) ,  a\acute{a} \in E(\Gamma)}(d_{\Gamma \times \Omega}(a,b)+d_{\Gamma \times \Omega}(\acute{a},b) - 2)^2 $\\

$\hspace*{2.2cm}$ $+\sum_{(\acute{a},\acute{b})(\acute{a},\acute{\acute{b}}) ,  \acute{b}\acute{\acute{b}} \in E(\Omega)}(d_{\Gamma \times \Omega}(\acute{a},\acute{b})+d_{\Gamma \times \Omega}(\acute{a},\acute{\acute{b}}) - 2)^2 $\\

$\hspace*{2.2cm}$ $+\sum_{(\acute{a},\acute{b})(\acute{\acute{a}},\acute{b}) , \acute{a}\acute{\acute{a}}\in E(\Gamma)}(d_{\Gamma \times \Omega}(\acute{a},\acute{b})+d_{\Gamma \times \Omega}(\acute{\acute{a}},\acute{b}) - 2)^2 $\\

$\hspace*{2.2cm}$ =$A_1+A_2+A_3+A_4$.\\

We consider\\

$\hspace*{1.6cm}$ $A_1=\sum_{(a,b)(a,\acute{b}) ,  b\acute{b} \in E(\Omega)}(d_{\Gamma \times \Omega}(a,b)+d_{\Gamma \times \Omega}(a,\acute{b}) - 2)^2 $\\

$\hspace*{2.2cm}$ =$\sum_{a\in V(\Gamma)}\sum_{b\acute{b}\in E(\Omega)}(2d_\Gamma (a)+d_\Omega (b)+d_\Omega (\acute{b}) - 2)^2$\\

$\hspace*{2.2cm}$ $=4\sum_{b\acute{b} \in E(\Omega)} \sum_ {a\in V(\Gamma)}(d_\Gamma(a)^2)$\\ 

$\hspace*{2.2cm}$ $+ \sum_{a\in V(\Gamma)} \sum_{ b\acute{b} \in E(\Omega)}(d_\Omega(b) + d_\Omega(\acute{b}) - 2)^2$\\

$\hspace*{2.2cm}$ $+4 \sum_{a \in V(\Gamma)} d_\Gamma (a) \sum_{b \acute{b} \in E(\Omega)}(d_\Omega (b)+d_\Omega (\acute{b}) - 2)$\\

$\hspace*{2.2cm}$ $=n_\Gamma EM_1(\Omega)+8m_\Gamma M_1(\Omega)+4m_\Omega M_1(\Gamma) - 16m_\Gamma m_\Omega.$\\

Similarly,\\

$\hspace*{1.6cm}$ $A_2=\sum_{(a,b)(\acute{a},b) , a\acute{a} \in E(\Gamma)}(d_{\Gamma \times \Omega}(a,b)+d_{\Gamma \times \Omega}(\acute{a},b) - 2)^2 $\\

$\hspace*{2.2cm}$ =$\sum_{b\in V(\Omega)}\sum_{a\acute{a}\in E(\Gamma)}(2d_\Omega (b)+d_\Gamma (a)+d_\Gamma (\acute{a}) - 2)^2$\\

$\hspace*{2.2cm}$ $=4\sum_{a\acute{a} \in E(\Gamma)}\sum_{ b\in V(\Omega)}(d_\Omega (b)^2)$\\ 

$\hspace*{2.2cm}$ $+ \sum_{b\in V(\Omega)} \sum_{ a\acute{a} \in E(\Gamma)}(d_\Gamma (a) + d_\Gamma (\acute{a}) - 2)^2$\\

$\hspace*{2.2cm}$ $+4 \sum_{b \in V(\Omega)} d_\Omega (b) \sum_{a \acute{a} \in E(\Gamma)}(d_\Gamma (a)+d_\Gamma (\acute{a}) - 2)$\\

$\hspace*{2.2cm}$ $=n_\Omega EM_1(\Gamma)+8m_\Omega M_1(\Gamma)+4m_\Gamma M_1(\Omega) - 16m_\Gamma m_\Omega.$\\

Similarly again,\\

$\hspace*{1.6cm}$ $A_3=\sum_{(\acute{a},\acute{b})(\acute{a},\acute{\acute{b}}) ,  \acute{b}\acute{\acute{b}} \in E(\Omega)}(d_{\Gamma \times \Omega}(\acute{a},\acute{b})+d_{\Gamma \times \Omega}(\acute{a},\acute{\acute{b}}) - 2)^2 $\\

$\hspace*{2.2cm}$ =$\sum_{\acute{a}\in V(\Gamma)}\sum_{\acute{b}\acute{\acute{b}}\in E(\Omega)}(2d_\Gamma (\acute{a})+d_\Omega (\acute{b})+d_\Omega (\acute{\acute{b}}) - 2)^2$\\

$\hspace*{2.2cm}$ $=4\sum_{\acute{b}\acute{\acute{b}} \in E(\Omega)} \sum_ {\acute{a}\in V(\Gamma)}(d_\Gamma(\acute{a})^2)$\\ 

$\hspace*{2.2cm}$ $+ \sum_{\acute{a}\in V(\Gamma)} \sum_{\acute{b}\acute{\acute{b}} \in E(\Omega)}(d_\Omega(\acute{b}) + d_\Omega(\acute{\acute{b}}) - 2)^2$\\

$\hspace*{2.2cm}$ $+4 \sum_{\acute{a} \in V(\Gamma)} d_\Gamma (\acute{a}) \sum_{\acute{b}\acute{\acute{b}} \in E(\Omega)}(d_\Omega (\acute{b})+d_\Omega (\acute{\acute{b}}) - 2)$\\

$\hspace*{2.2cm}$ $=n_\Gamma EM_1(\Omega)+8m_\Gamma M_1(\Omega)+4m_\Omega M_1(\Gamma) - 16m_\Gamma m_\Omega.$\\

Finally,\\

$\hspace*{1.6cm}$ $A_4=\sum_{(\acute{a},\acute{b})(\acute{\acute{a}},\acute{b}) , \acute{a}\acute{\acute{a}} \in E(\Gamma)}(d_{\Gamma \times \Omega}(\acute{a},\acute{b})+d_{\Gamma \times \Omega}(\acute{\acute{a}},\acute{b}) - 2)^2 $\\

$\hspace*{2.2cm}$ =$\sum_{\acute{b}\in V(\Omega)}\sum_{\acute{a}\acute{\acute{a}}\in E(\Gamma)}(2d_\Omega (\acute{b})+d_\Gamma (\acute{a})+d_\Gamma (\acute{\acute{a}}) - 2)^2$\\

$\hspace*{2.2cm}$ $=4\sum_{\acute{a}\acute{\acute{a}} \in E(\Gamma)}\sum_{ \acute{b}\in V(\Omega)}(d_\Omega (\acute{b})^2)$\\ 

$\hspace*{2.2cm}$ $+ \sum_{\acute{b}\in V(\Omega)} \sum_{\acute{a}\acute{\acute{a}} \in E(\Gamma)}(d_\Gamma (\acute{a}) + d_\Gamma (\acute{\acute{a}}) - 2)^2$\\

$\hspace*{2.2cm}$ $+4 \sum_{\acute{b} \in V(\Omega)} d_\Omega (\acute{b}) \sum_{\acute{a}\acute{\acute{a}} \in E(\Gamma)}(d_\Gamma (\acute{a})+d_\Gamma (\acute{\acute{a}}) - 2)$\\

$\hspace*{2.2cm}$ $=n_\Omega EM_1(\Gamma)+8m_\Omega M_1(\Gamma)+4m_\Gamma M_1(\Omega) - 16m_\Gamma m_\Omega.$\\

By adding $A_1,A_2,A_3$ and $+A_4$, the desired result follows.

\end{proof}

\section{Methods and Theory}
In this section, we are going to compute the edge hyper-Zagreb index for the family of linear Acenes $(C_{4n+2}H_{2n+4})$ molecule, (FIGURE 1).

\begin{figure}[H]
\begin{center}
\includegraphics[scale=0.46]{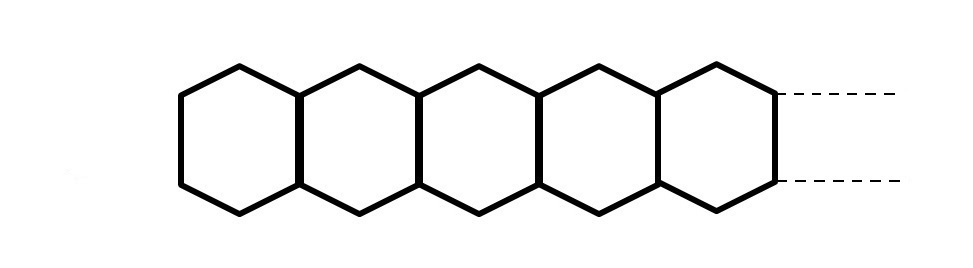}
\caption{$Linear$ $Acenes$ $Molecular$ $Graph$}
\label{fig:1} 
\end{center}      
\end{figure}

\begin{theorem}
If $\Gamma$ is a linear Acenes $(C_{4n+2}H_{2n+4})$ molecule, then its edge hyper-Zagreb index is equal to
\begin{center}
$EHM(\Gamma)=4(85n-62).$
\end{center}
\end{theorem} 
 
\begin{proof}
Let $\Gamma$ be a simple graph of Acenes family. The rings set divide into two partitions based on outer (pendant) rings and inner rings. Now, we redivide rings according to the degree of edges as follow (see FIGURE 2 and FIGURE 3):\\
$(i)$ degree of edges in outer rings. So, we have 
\begin{center}
$\acute{E}_1=\{\alpha\thicksim \beta \in E(\Gamma)|d_\alpha=d_\beta=2\}$,
\end{center}
\begin{center}
$\acute{E}_2=\{\alpha\thicksim \beta \in E(\Gamma)|d_\alpha=2, d_\beta=3\}$ and
\end{center}
\begin{center}
$\acute{E}_3=\{\alpha\thicksim \beta \in E(\Gamma)|d_\alpha=3, d_\beta=4\}$.
\end{center}
$(ii)$ degree og edges in inner rings. So, we obtain
\begin{center}
$\acute{E}_4=\{\alpha\thicksim \beta \in E(\Gamma)|d_\alpha=d_\beta=3\}$ and
\end{center}
\begin{center}
$\acute{E}_5=\{\alpha\thicksim \beta \in E(\Gamma)|d_\alpha=3, d_\beta=4\}$.
\end{center}
$(iii)$ degree of edges between outer and inner rings. Then, we attain
\begin{center}
$\acute{E}_6=\{\alpha\thicksim \beta \in E(\Gamma)|d_\alpha= d_\beta=3\}$.
\end{center}

\begin{figure}[H]
\begin{center}
\includegraphics[scale=0.43]{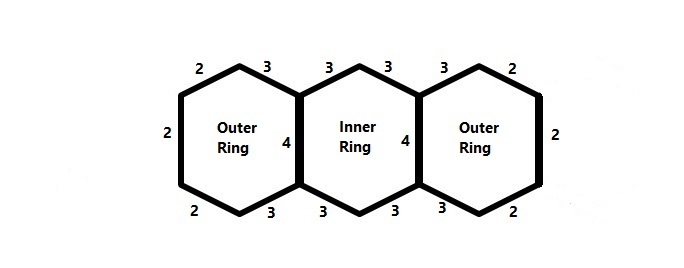}
\caption{$Phenantherene$ $Molecular$ $Graph$}
\label{fig:2}
\end{center}       
\end{figure}

\begin{figure}[H]
\begin{center}
\includegraphics[scale=0.43]{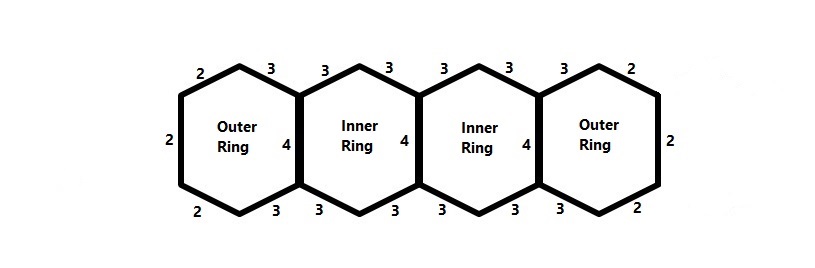}
\caption{$Chrysene$ $Molecular$ $Graph$}
\label{fig:3}
\end{center}       
\end{figure}

Now, we investigate the statement by induction on the number of rings that it is indicated by $n$. According to FIGURE 1 and FIGURE 2, we can compute $|\acute{E}_1|,...,|\acute{E}_6|$ for $n=3$ and $n=4$.\\
If $n=3$, then 
\begin{center}
$|\acute{E}_1(\Gamma)|=4, |\acute{E}_2(\Gamma)|=4, |\acute{E}_3(\Gamma)|=4, |\acute{E}_4(\Gamma)|=2, |\acute{E}_5(\Gamma)|=4$ and $|\acute{E}_6(\Gamma)|=4$.
\end{center}
So, we can conclude $EHM(\Gamma)$ as follow:
\begin{center}
$EHM(\Gamma)= \sum_{\alpha\thicksim \beta \in \acute{E}_i(\Gamma)} (d_\alpha+d_\beta)^2$, for $1\leqslant i\leqslant 6\Rightarrow$
\end{center}
\begin{center}
$EHM(\Gamma)= \sum_{\alpha\thicksim \beta \in \acute{E}_1(\Gamma)} (d_\alpha+d_\beta)^2 + \sum_{\alpha\thicksim \beta \in \acute{E}_3(\Gamma)} (d_\alpha+d_\beta)^2 + \sum_{\alpha\thicksim \beta \in \acute{E}_4(\Gamma)} (d_\alpha+d_\beta)^2 + \sum_{\alpha\thicksim \beta \in \acute{E}_5(\Gamma)} (d_\alpha+d_\beta)^2 + \sum_{\alpha\thicksim \beta \in \acute{E}_6(\Gamma)} (d_\alpha+d_\beta)^2$
\end{center}
\begin{center}
$= 4(2+2)^2+4(2+3)^2+4(3+4)^2+2(3+3)^2+4(3+4)^2+4(3+3)^2=772$.
\end{center}
If $n=4$, then 
\begin{center}
$|\acute{E}_1(\Gamma)|=4, |\acute{E}_2(\Gamma)|=4, |\acute{E}_3(\Gamma)|=4, |\acute{E}_4(\Gamma)|=4, |\acute{E}_5(\Gamma)|=8$ and $|\acute{E}_6(\Gamma)|=6$.
\end{center}
Furthermore, same as previous calculations $EHM(\Gamma)$ for $n=3$, we can verify $EHM(\Gamma)$ for $n=4$. Hence, $EHM(\Gamma)=1112$.\\

Now, let $n>4$.\\
We assume that the claim holds for the number of rings less than $n$. We saw that $|\acute{E}_1(\Gamma)| = |\acute{E}_2(\Gamma)| = |\acute{E}_3(\Gamma)|=4$ are constant for all edges in outer rings. So,
\begin{center} 
$EHM(\Gamma)= \sum_{\alpha\thicksim \beta \in \acute{E}_i(\Gamma)} (d_\alpha+d_\beta)^2 = 360$ for $1\leqslant i\leqslant 3$
\end{center}
and 
\begin{center}
$|\acute{E}_4(\Gamma)|=2(n-2), |\acute{E}_5(\Gamma)|=4(n-2)$ and $|\acute{E}_6(\Gamma)|=2(n-1)$, 
\end {center}
in which $n$ is the number of rings. Hence, the edge hyper-Zagreb index of Acenes $(C_{4n+2}H_{2n+4})$ is equal to 
\begin{center}
$EHM(\Gamma)= \sum_{\alpha\thicksim \beta \in \acute{E}_i(\Gamma)} (d_\alpha+d_\beta)^2$, for $1\leqslant i\leqslant 6\Rightarrow$
\end{center}
\begin{center}
$EHM(\Gamma)= 360 + \sum_{\alpha\thicksim \beta \in \acute{E}_4(\Gamma)} (d_\alpha+d_\beta)^2 + \sum_{\alpha\thicksim \beta \in \acute{E}_5(\Gamma)} (d_\alpha+d_\beta)^2 + \sum_{\alpha\thicksim \beta \in \acute{E}_6(\Gamma)} (d_\alpha+d_\beta)^2$
\end{center}
\begin{center}
$= 360+2(n-2)(6)^2+4(n-2)(7)^2+2(n-1)(6)^2$
\end{center}
\begin{center}
$= 4(85n-62).$
\end{center}
Finally, the proof is completed.
 \end{proof}
 
 \section{Results}
 
Equation $EHM(A)=4(85n-62)$ was used to calculate the edge hyper-Zagreb index for the eight elements of the Acenes family, which are then displayed in TABLE 1.\\
\begin{table}[H]
\centering
\caption{$EHM$ index for the first eight members of Acenes family}
\begin{tabular}{|c|c|c|}
\hline
Chemical Formula & IUPAC Name & $EHM(\Gamma)$\\
\hline
$C_{10}H_8$ & Benzene & $432$ \\
\hline
$C_{14}H_{10}$ & Naphthalene & $772$ \\
\hline
$C_{18}H_{12}$ & Anthracene & $1112$\\
\hline
$C_{22}H_{14}$ & Tetracene & $1452$\\
\hline
$C_{26}H_{16}$ & Pentacene & $1792$\\
\hline
$C_{30}H_{18}$ & Hexacene & $2132$\\
\hline
$C_{34}H_{20}$ & Heptacene & $2472$\\
\hline
$C_{36}H_{22}$ & Optacene & $2812$\\
\hline
\end{tabular}
\end{table}
We conducted calculations for Gibbs Energy $(GE)$, Heat of formation $(H.oF)$, as well as Electro-optical characteristics like gap energy $(E_g)$ and electron affinity energy $(E_{ea})$ in the Acenes family using empirical data to examine the thermodynamic properties. These results were then compared against established sources, as detailed in TABLE 2 [10,12].
\begin{table}[H]
\centering
\caption{$(GE)$, $(H.oF)$, $(E_g)$ and $(E_{ea})$ energy}
\begin{tabular}{|c|c|c|c|c|}
\hline
Chemical Formula & $H.of (Kj/mol)$ & $GE (Kj/mol)$ & $E_g$ & $E_{ea}$\\
\hline
$C_{10}H_8$ & $80.83$ & $121.68$ & $-0.32$ & $4.09$\\
\hline
$C_{14}H_{10}$ & $177.87$ & $252.38$ & $-0.05$ & $4.19$\\
\hline
$C_{18}H_{12}$ & $274.91$ & $383.08$ & $-0.29$ & $3.73$\\
\hline
$C_{22}H_{14}$ & $371.95$ & $513.78$ & $0.4$ & $3.7$\\
\hline
$C_{26}H_{16}$ & $468.99$ & $644.48$ & $0.57$ & $3.47$\\
\hline
$C_{30}H_{18}$ & $566.03$ & $775.18$ & $ 0.64$ & $3.5$\\
\hline
$C_{34}H_{20}$ & $633.07$ & $905.88$ & $0.73$ & $3.44$\\
\hline
$C_{36}H_{22}$ & $760.11$ & $1036.58$ & $ 0.84$ & $3.36$\\
\hline
\end{tabular}
\end{table}
FIGURES 4, 5, 6 and 7 illustrate the values of $E_{ea}$, $E_g$, $H.oF$, and $GE$ for linear Acenes in relation to the edge hyper-Zagreb index, reflecting the findings of this study.\\ 
FIGURES 4 and 5 demonstrate the precise calculation of $E_g$ and $E_{ea}$ using the $EHM(A)$ index variations within the Acenes family.

\begin{figure}[H]
\centering
\includegraphics[scale=0.48]{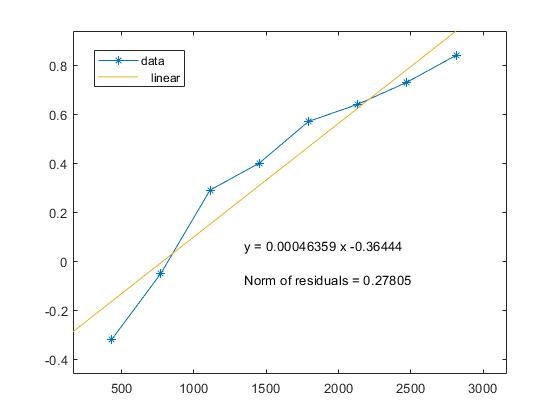}
\caption{Electron affinity energy}
\label{fig:1}       
\end{figure}
\begin{figure}[H]
\centering
\includegraphics[scale=0.48]{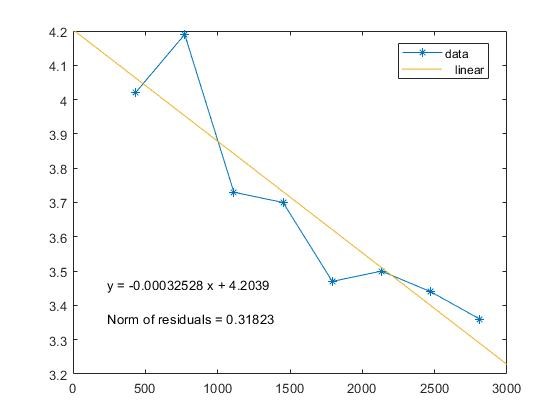}
\caption{Gap energy}
\label{fig:1}       
\end{figure}
\begin{figure}[H]
\centering
\includegraphics[scale=0.48]{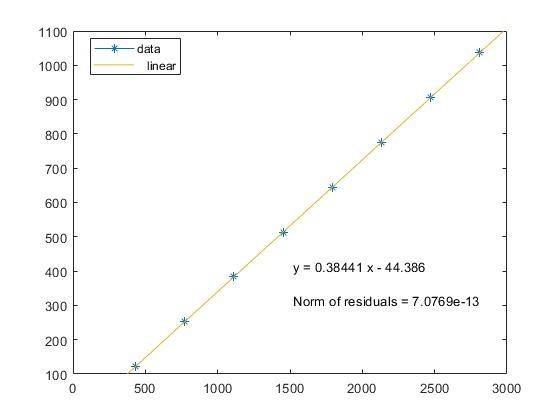}
\caption{Gipps energy}
\label{fig:1}       
\end{figure}
\begin{figure}[H]
\centering
\includegraphics[scale=0.48]{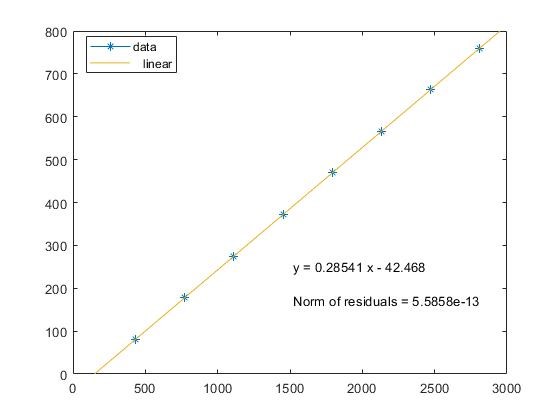}
\caption{Heat of formation energy}
\label{fig:1}       
\end{figure}
These forecasts can be made using a pair of equations: 
\begin{center}
$E_{ea}=0.00046359 EHM(\Gamma) - 0.36444$ $\hspace*{2cm}(4-1)$,       
\end{center}
\begin{center}
$E_g=-0.00032528 EHM(\Gamma) +4.2039$ $\hspace*{2cm}(4-2)$. 
\end{center}
FIGURES (6) and (7) demonstrate the effectiveness of the $EHM(\Gamma)$ index in predicting $GE$ and $H.oF$ within the linear Acenes family, with an $R^2$ value of $1$. The corresponding equation is provided below 
\begin{center}
$GE=0.38441 EHM(\Gamma) - 44.386$ $\hspace*{2.7cm}(4-3)$,       
\end{center}
\begin{center}
$E_g=-0.00032528 EHM(\Gamma) +4.2039$ $\hspace*{2cm}(4-4)$.      
\end{center}
\section{Conclusion}
In TABLE 3, we present the values of $E_G$, $E_{ea}$, $GE$, and $H.oF$ for the Acenes family using the $TIM$ method. This table highlights the reliability of the $TIM$ method by comparing the results with the reference values from TABLE 2. The $TIM$ method, along with the edge hyper-Zagreb index, allows for accurate prediction of various physical and chemical properties within the Acenes family, saving time and resources compared to traditional theoretical and experimental approaches that often yield only approximate results. Moving forward, several heavier members of the Acenes family will undergo analysis using the $TIM$ method. Equations $(4-1), (4-2), (4-3)$ and $(4-4)$ have been employed to forecast the energies of $E_G$, $E_{ea}$, $GE$, and $H.oF$ with the outcomes detailed in TABLE 4.
\begin{table}[H]
\centering
\caption{Calculation $(GE)$, $(H.oF)$, $(E_g)$ and $(E_{ea})$ energy of first eight Acenes family by $TIM$ method}
\begin{tabular}{|c|c|c|c|c|}
\hline
Chemical Formula & $H.of$ & $GE$ & $E_g$ & $E_{ea}$\\
\hline
$C_{10}H_8$ & $80.82912$ & $121.67912$ & $-0.16412912$ & $4.06337904$\\
\hline
$C_{14}H_{10}$ & $177.86852$ & $252.37852$ & $-0.00650852$ & $3.95278384$\\
\hline
$C_{18}H_{12}$ & $274.90792$ & $383.07792$ & $0.15111208$ & $3.84218864$\\
\hline
$C_{22}H_{14}$ & $371.94732$ & $513.77732$ & $0.30873268$ & $3.73159344$\\
\hline
$C_{26}H_{16}$ & $468.98672$ & $644.47672$ & $0.46635328$ & $3.62099824$\\
\hline
$C_{30}H_{18}$ & $566.02612$ & $775.17612$ & $ 0.62397388$ & $3.51040304$\\
\hline
$C_{34}H_{20}$ & $633.06552$ & $905.87552$ & $0.78159448$ & $3.39980784$\\
\hline
$C_{36}H_{22}$ & $760.10492$ & $1036.57492$ & $ 0.93921508$ & $3.28921264$\\
\hline
\end{tabular}
\end{table}
 
\begin{table}[H]
\centering
\caption{Prediction $(GE)$, $(H.oF)$, $(E_g)$ and $(E_{ea})$ energy for Acenes family}
\begin{tabular}{|c|c|c|c|c|c|c|}
\hline
Chemical Formula & $n$ & $EHM(\Gamma)$ & $H.oF$ & $GE$ & $E_{ea}$ & $E_g$\\
\hline
$C_{42}H_{24}$ & $10$ & $3152$ & $80.82912$ & $121.67912$ & $-0.16412912$ & $4.06337904$\\
\hline
$C_{46}H_{26}$ & $11$ & $3492$ & $177.86852$ & $252.37852$ & $-0.00650852$ & $3.95278384$\\
\hline
$C_{50}H_{28}$ & $12$ & $3832$ & $274.90792$ & $383.07792$ & $0.15111208$ & $3.84218864$\\
\hline
$C_{54}H_{30}$ & $13$ & $4172$ & $371.94732$ & $513.77732$ & $0.30873268$ & $3.73159344$\\
\hline
$C_{58}H_{32}$ & $14$ & $4512$ & $468.98672$ & $644.47672$ & $0.46635328$ & $3.62099824$\\
\hline
$C_{62}H_{34}$ & $15$ & $4852$ & $566.02612$ & $775.17612$ & $0.62397388$ & $3.51040304$\\
\hline
$C_{66}H_{36}$ & $16$ & $5192$ & $663.06552$ & $905.87552$ & $0.78159448$ & $3.39980784$\\
\hline
$C_{70}H_{38}$ & $17$ & $5532$ & $760.10492$ & $1036.57492$ & $0.93921508$ & $3.28921264$\\
\hline
\end{tabular}
\end{table}

 \noindent Zohreh Aliannejadi \\ Department of Mathematics, Islamic Azad University, South Tehran Branch, \\ Tehran, Iran. \\ e-mail: z\_alian@azad.ac.ir \\ phone: +989122797710
 \bigskip \noindent 
 \\Somayeh Shafiee Alamoti \\ Department of Mathematics, Islamic Azad University, South Tehran Branch, \\ Tehran, Iran. \\ e-mail: shafiee.s88@gmail.com \\ Phone: +989126994492
 \bigskip \noindent 

 \end{document}